\documentclass[11pt,a4paper]{article}
\usepackage{amsmath}
\usepackage{amssymb}
\usepackage{bm}

\renewcommand{\Re}{\mathop{\rm Re}\nolimits}

\newtheorem{definition}{Definition}[section]
\newtheorem{thm}{Theorem}[section]

\newtheorem{lem}[thm]{Lemma}

\newtheorem{rem}{Remark}[section]

\newcommand{\intd}{\textrm{d}}

\newcommand{\dps}{\displaystyle}

\newcommand{\eval}[2][\right]{\relax
  \ifx#1\right\relax \left.\fi#2#1\rvert}

\begin{document}

\title{Global solutions for a nonlocal Ginzberg-Landau equation and a nonlocal Fokker-Plank equation
\footnote{Supported partly by the   NSF Grant 1025422, the NSFC grant 11271290 and 11171158 and   National Basic Research Program of China (973
Program) No. 2013CB834100.} }
\author{Jinchun He\thanks{E-mail address: taoismnature@mail.hust.edu.cn (J. He)}, Jinqiao Duan\thanks{E-mail address: duan@iit.edu (J. Duan)} and  Hongjun Gao\thanks{E-mail address: gaohj@njnu.edu.cn (H. Gao, corresponding author)}}
\date{December 8, 2013}
\maketitle

\begin{center}
$\ast$~School of Mathematics and Statistics, Huazhong University of Science and Technology, Wuhan 430074, China. \\
$\dag$~Department of Applied Mathematics, Illinois Institute of Technology, Chicago, IL 60616, USA. \\
$\ddag$~Institute of Mathematics, Nanjing Normal University, Nanjing 210023, China \\ and Institute of Mathematics, Jilin University, Changchun 130012, China
\end{center}

\textbf{Abstract:}~~This work is devoted to the study of a nonlocal Ginzberg-Landau equation by the semigroup method and a nonlocal Fokker-Plank equation by the viscosity vanishing method. For the nonlocal Ginzberg-Landau equation, there exists  a unique global solution in the set $C^0(\mathbb{R}^+,\,H_0^{\frac{\alpha}{2}}(D))\cap L_{loc}(\mathbb{R}^+,\,H_0^{\alpha}(D))$, for $\alpha\in (0,\,2)$. For the nonlocal Fokker-Plank equation, the regularity of the solution is weaker than that of the nonlocal Ginzberg-Landau equation due to the drift term.


\section{Introduction}

For a system described by a scalar stochastic differential equation with Brownian motion (a Gaussian process),
$$ dX_t =b(X_t)dt + dB_t, \;\; X_0 =x$$
the corresponding Fokker-Planck equation contains the usual Laplacian operator $\Delta$,
$$
p_t = \frac12 \Delta p -(b(x)p)_x.
$$
If the Brownian motion is replaced by a $\alpha$-stable L\'evy motion (a non-Gaussian process) $L_t^\alpha$ for $\alpha \in (0, 2)$, the Fokker-Planck equation becomes a nonlocal partial differential equation \cite{Applebaum} with a fractional Laplacian operator $(-\Delta)^{\frac{\alpha}{2}}$. When the drift term  $b$ in the   stochastic differential equation depends on the probability distribution of the system state, the Fokker-Planck equation becomes nonlinear.

Nonlocal partial differential equations also arise in the modeling of anomalous diffusion of particles in fluids \cite{MS2012, MK2004}. Moreover, nonlocal effects and long range diffusion occur in the cell density evolution in certain biological processes, such as embryological development \cite[Section 11.5]{Murray}.

It is desirable to investigate analytical foundation of this type of nonlocal partial differential equations.
In the present paper, we consider the following time-dependent nonlocal
    partial differential equation
\begin{equation}\label{mysystem}
\left\{\begin{array}{ll}
u_t=A_{\alpha}u+F(u,\,\nabla u),& ~~x\in D,~t>0 \\
u|_{D^c}=0, &\\
u(x,0)=u_0(x), &
\end{array}
\right.
\end{equation}
where $D=(0,1)$ is an interval in $\mathbb{R}^1$, $D^c=\mathbb{R}^1\setminus D$, and $A_{\alpha}=-(-\Delta)^{\alpha/2}$ is the nonlocal Laplacian operator with definition
$$
(-\Delta)^{\alpha/2} f(x) = c_{\alpha} \int_{\mathbb{R}^1\setminus\{0\}} \frac{f(x+y)-f(x)} {|y|^{1+\alpha}} \intd y,~~0<\alpha<2,
$$
where $c_{\alpha}$ is a constant depending on the order $\alpha$. Volume constraints are natural extensions to the nonlocal case of boundary conditions for differential equations.

This paper is organized as follows. In Section 2, we study the local solution   for the nonlocal Ginzberg-Landau equation in the space $C^0([0,\,T],\,H_0^{\alpha/2}(D))$, based on several properties of the semigroup generated by the nonlocal Laplacian operator $A_{\alpha}$. We also obtain a further result  for the global solution by the energy estimates. Moreover, the same method   also be applies to   other cases where $u-u^3$ is replaced by a dissipative term $F(u)$.
Section 3 focuses on a nonlocal diffusion equation with a drift term. Equations with drift and anomalous diffusion appear in numerous places in mathematical physics.
A successful understanding of well-posedness of the problem relies on the a priori estimates that will be established. In many  cases, these are based on the linearized drift-diffusion equation, and this  provides the motivation for the present work. Hence, a nonlocal Fokker-Plank equation is considered in this paper. Comparing with Caffarelli and Silvestre's papers \cite{SVZ2013}--\cite{CS2007} with spatial regularity specified  in Banach spaces, we work in spaces where spatial regularity is in   Hilbert spaces. In order to overcome the difficulties originating from the drift term, the viscosity vanishing method is applied.

\section{A nonlocal Ginzberg-Landau equation}

We consider the following nonlocal Ginzberg-Landau equation by semigroup method \cite{DH1981, Pazy1983}
\begin{equation}\label{semilin}
\left\{\begin{array}{ll}
u_t=A_{\alpha}u+u-u^3,& ~~x\in D=(0, 1), ~t>0 \\
u|_{D^c}=0, &\\
u(x,0)=u_0(x). &
\end{array}
\right.
\end{equation}
In order to solve the problem ($\ref{semilin}$), we need the properties of the semigroup generated by the nonlocal Laplacian operator $A_{\alpha}$ and some formulas of nonlocal calculus \cite{DGLZ2011}.

\subsection{Some estimates on the nonlocal Laplacian}\label{semigroup-esti}

\begin{definition}
$A$ is called a sectorial operator if (i) it is dense defined, (ii) for some $\phi\in(0,\,\frac{\pi}{2})$,
$M\ge 1$, and $a\in\mathbb{R}$, $S_{a,\,\phi}=\{\lambda|~\phi\le|\mbox{\textrm{arg}}(\lambda-a)|\le\pi,~\lambda\neq a\}
\subset\rho(A)$, with resolvat set $\rho(A)$ and (iii) $\|(\lambda I-A)^{-1}\|\le M/|\lambda-a|$.
\end{definition}

\begin{lem}(\cite{Kw2012})
Let $(-\Delta)^{\frac{\alpha}{2}}$ be defined in $L^2(D)$ for $\alpha \in (0, 2)$. Then the eigenvalues of the spectral problem
\begin{equation}
(-\Delta)^{\frac{\alpha}{2}}\varphi(x)=\lambda\varphi(x),~x\in  D,
\end{equation}
where $\varphi(x)\in L^2(D)$ is extended to $\mathbb{R}$ by $0$ is
\begin{equation}
\lambda_n=(\frac{n}{2}-\frac{(2-\alpha)}{8})^\alpha+o(\frac{1}{n}),
\end{equation}
and $\lambda_n$ satisfies
$$0 < \lambda_1 < \lambda_2 \leq \cdot\cdot\cdot\leq \lambda_j
\leq \cdot\cdot\cdot,$$
and the corresponding eigenfunctions $\varphi_n$, after an appropriate normalization, form a complete orthonormal basis  in $L^2(D)$.
\end{lem}

Using the above lemmas, it can be verified that the nonlocal Laplacian operator $A_{\alpha}$ is sectorial,  and the following estimates hold.
\begin{lem}\label{semiest}
The nonlocal Laplacian operator $A_{\alpha}$ satisfies the estimates
\begin{equation}
\|\textrm{e}^{tA_{\alpha}}\|_{L^2(D)}\le C\textrm{e}^{-\delta t},~\mbox{and}~~~\|A_{\alpha}\textrm{e}^{tA_{\alpha}}\|_{L^2(D)}
\le \frac{C}{t}\textrm{e}^{-\delta t},
\end{equation}
where $C,\,\delta $ are positive constants independent of $t$.
\end{lem}
\textbf{Proof}~~Set $\mu=\lambda t$ for $\lambda>0$. Then
$$
\|\textrm{e}^{tA_{\alpha}}\|_{L^2(D)}=\dps{\left\|\frac{1}{2\pi i}\int_{\Gamma}\textrm{e}^{\mu}(\frac{\mu}{t}-A_{\alpha})^{-1}
\frac{\intd\mu}{t}\right\|}_{L^2(D)} \le \frac{M}{2\pi}\int_{\Gamma}\mid\textrm{e}^{\mu}\mid
\frac{\mid\intd\mu\mid}{\mid \mu\mid} \le C\textrm{e}^{-\delta t}.
$$
Moreover,
$$
\|A_{\alpha}\textrm{e}^{tA_{\alpha}}\|_{L^2(D)}=\dps{\left\|\frac{1}{2\pi i}A_{\alpha}\int_{\Gamma}\textrm{e}^{\mu}(\frac{\mu}{t}-A_{\alpha})^{-1}
\frac{\intd\mu}{t}\right\|}_{L^2(D)} \le \frac{1}{2\pi}\frac{M}{\delta}\int_{\Gamma}\mid\textrm{e}^{\mu}\mid
\frac{\mid\intd\mu\mid}{\mid \mu\mid}\frac{1}{t} \le \frac{C}{t}\textrm{e}^{-\delta t}.
$$
This proves the lemma. $\Box$

\begin{definition}\label{indexoperator}
Let $A$ be a sectorial operator and $\Re \sigma(-A)>0$.   For every $\beta>0$, define
$A^{-\beta}=\dps{\frac{1}{\Gamma(\beta)}\int_0^{\infty}t^{\beta-1}\textrm{e}^{tA}\intd t}$.
Moreover, $A^{\beta}=(A^{-\beta})^{-1}$, and $A^0=I$.
\end{definition}

\begin{lem}\label{semiregu}
Let $A$ be a sectorial operator, and $\Re \sigma(-A)>\delta>0$. Then $\forall\, \beta\ge 0$,
$\exists\,C(\beta)<\infty$, such that $\forall\,t>0$,
\begin{equation}
\|A^{\beta}\textrm{e}^{tA}\|\le C(\beta) t^{-\beta}\textrm{e}^{-\delta t}.
\end{equation}
\end{lem}
\textbf{Proof}~~$\forall$ $m=1,\,2,\,\cdots$,
\begin{equation}
\|A^m\textrm{e}^{tA}\|=\left\|\left(A\textrm{e}^{\frac{tA}{m}}\right)^m\right\|
\le (Cm)^m t^{-m}\textrm{e}^{-\delta t}.
\end{equation}
For $0<\beta<1$ and  $t>0$,
\begin{eqnarray*}
\|A^{\beta}\textrm{e}^{tA}\|&=&\|A^{\beta-1}A\textrm{e}^{tA}\|=\|A^{-(1-\beta)}A\textrm{e}^{tA}\| \\
&=&\dps{\left\|\frac{1}{\Gamma(1-\beta)}\int_0^{\infty}\tau^{1-\beta-1}A\textrm{e}^{-A(t+\tau)}
\intd \tau\right\|}\le \frac{1}{\Gamma(1-\beta)}\int_0^{\infty}\tau^{-\beta}
\|A\textrm{e}^{-A(t+\tau)}\| \intd \tau \\
&\le& \frac{C}{\Gamma(1-\beta)}\int_0^{\infty}\tau^{-\beta}(t+\tau)^{-1}\textrm{e}^{-\delta(t+\tau)}
\intd \tau=C\Gamma(\beta)t^{-\beta}\textrm{e}^{-\delta t}.
\end{eqnarray*}
Hence, $\forall$ $\beta\ge 0$, $\|A^{\beta}\textrm{e}^{tA}\|\le C(\beta) t^{-\beta}\textrm{e}^{-\delta t}$.
$\Box$

\begin{lem}\label{embetoHs}
$\textrm{Dom}(A_{\alpha}) \hookrightarrow H_0^{\frac{\alpha}{2}}(D)$.
\end{lem}
\textbf{Proof}~~From \cite{CMN2012}, we know that $\textrm{Dom}(A_{\alpha}) \subset H_0^{\frac{\alpha}{2}}(D)$. According to the embedding
results and the nonlocal calculus in   \cite{DGLZ2011}, and using the H\"{o}lder inequality,  we conclude that
\begin{eqnarray*}
\|u\|^2_{H^{\frac{\alpha}{2}}(D)}&=&\|u\|^2_{L^2(D)}+|u|^2_{H^{\frac{\alpha}{2}}(D)} \\
&\le& \|u\|^2_{L^2(D)}+\frac12\int_D\int_D(\mathcal{D}^{\ast}(u)(x,\,y))^2 \intd x \intd y
+C\varepsilon^{-(1+\alpha)}\|u\|^2_{L^2(D)}\\
&=& C(D,\,\varepsilon,\,s)\|u\|^2_{L^2(D)}-\frac12 \langle A_{\alpha}u,\,u \rangle \le C\|u\|^2_{L^2(D)}
+C\|A_{\alpha}u\|_{L^2(D)}\|u\|_{L^2(D)} \\
&=& C\|u\|_{L^2(D)}(\|u\|_{L^2(D)}+\|A_{\alpha}u\|_{L^2(D)}).
\end{eqnarray*}
By the nonlocal Poincar\`{e}'s inequality \cite{DGLZ2011}, we get
\begin{eqnarray*}
\|u\|^2_{H^{\frac{\alpha}{2}}(D)}&\le& C\|u\|_{H^{\frac{\alpha}{2}}(D)}(\|u\|_{L^2(D)}+\|A_{\alpha}u\|_{L^2(D)}).
\end{eqnarray*}
Therefore, $\|u\|_{H^{\frac{\alpha}{2}}(D)}\le C(\|u\|_{L^2(D)}+\|A_{\alpha}u\|_{L^2(D)})$. $\Box$

\begin{lem}\label{fracembe}
$\textrm{Dom}(A_{\alpha}^{\beta}) \hookrightarrow H^{\frac{\alpha}{2}}(D)$, $\frac12<\beta<1$.
\end{lem}
\textbf{Proof}~~By the Lemma $\ref{embetoHs}$, we have $\textrm{Dom}(A_{\alpha}^{\beta})
 \hookrightarrow \textrm{Dom}(A_{\alpha})  \hookrightarrow H^{\frac{\alpha}{2}}(D)$. $\Box$

\subsection{Local solution}

The existence   of the local solution comes from a standard contraction mapping argument. With number $T>0$ and $R>0$ to be fixed below, in the Banach space $X=C^0([0,\,T],\,H_0^{\frac{\alpha}{2}}(D))$, we consider a closed set
$$
S=\{u\in X:~\|u-u_0\|_X\le R\}.
$$
It follows that the mapping
\begin{equation}\label{solsemi}
u=\textrm{e}^{tA_{\alpha}}u_0+\int_0^t \textrm{e}^{(t-\tau)A_{\alpha}}(u-u^3)\intd \tau:=\Phi(u)
\end{equation}
is a contraction from $S$ into itself.

Since $\textrm{e}^{tA_{\alpha}}$ is a strongly continuous semigroup, we can choose $T_1$ such that
$\|\textrm{e}^{tA_{\alpha}}u_0-u_0\|_{H^{\frac{\alpha}{2}}(D)}\le R/2$ for $t\in[0,\,T_1]$. Now let $u\in X$.
Since $u^3$ is Lipschitz continuous from bounded subsets of $L^{6}(D)$ to $L^{2}(D)$,
we have a bound $\|u-u^3\|_X\le K_1 R$. Thus using Lemma ($\ref{semiest}$) and ($\ref{fracembe}$), we have
\begin{eqnarray}\label{intoS}
&&\int_0^t \|\textrm{e}^{(t-\tau)A_{\alpha}}(u-u^3)\|_{H^{\frac{\alpha}{2}}(D)}\intd \tau \nonumber \\
&\le& C\int_0^t \|\textrm{e}^{(t-\tau)A_{\alpha}}(u-u^3)\|_{L^2(D)}\intd \tau +C\int_0^t \|A^{\beta}_{\alpha}\textrm{e}^{(t-\tau)A_{\alpha}}(u-u^3)\|_{L^2(D)}\intd \tau \nonumber \\
&\le& C\int_0^t\|\textrm{e}^{(t-\tau)A_{\alpha}}\|_{L^2(D)}\|u-u^3\|_{L^2(D)}\intd \tau +C\int_0^t \|A^{\beta}_{\alpha}\textrm{e}^{(t-\tau)A_{\alpha}}\|_{L^2(D)}\|u-u^3\|_{L^2(D)}\intd \tau \nonumber \\
&\le& CK_1R\int_0^t\textrm{e}^{-\delta(t-\tau)}\intd \tau+CK_1R\int_0^t\frac{\textrm{e}^{-\delta(t-\tau)}}{(t-\tau)^{\beta}}\intd \tau \nonumber \\
&\le& \frac{C K_1R}{\delta}(1-\textrm{e}^{-\delta T_2})+\frac{C K_1R}{1-\beta}T_2^{1-\beta},
\end{eqnarray}
where $\frac12<\beta<1$. If we pick up $T_2\le T_1$ small enough, the right hand side of ($\ref{intoS}$) will be less than  $ R/2$ for $t\in[0,\,T_2]$. Therefore,
$\Phi:~S\rightarrow S$, provided $T\le T_2$.

To see  that $\Phi$ be a contraction, we also use the Lipschitz properties of $u-u^3$.
For $u,\,\bar{u}\in X$ and for $t\in [0,T_2]$, by Lemma ($\ref{semiest}$) and ($\ref{embetoHs}$), we have
\begin{eqnarray}
&&\|\Phi(u)-\Phi(\bar{u})\|_{H^{\frac{\alpha}{2}}(D)}\le\int_0^t \|\textrm{e}^{(t-\tau)A_{\alpha}}((u-\bar{u})-(u^3-\bar{u}^3))\|_{H^{\frac{\alpha}{2}}(D)}\intd \tau \nonumber \\
&\le& C\int_0^t \|\textrm{e}^{(t-\tau)A_{\alpha}}((u-\bar{u})-(u^3-\bar{u}^3))\|_{L^2(D)}\intd \tau  \nonumber \\
&&+C\int_0^t \|A^{\beta}_{\alpha}\textrm{e}^{(t-\tau)A_{\alpha}}((u-\bar{u})-(u^3-\bar{u}^3))\|_{L^2(D)}\intd \tau \nonumber \\
&\le& CL(R)\int_0^t\|\textrm{e}^{(t-\tau)A_{\alpha}}\|_{L^2(D)}\|u-\bar{u}\|_{L^2(D)}\intd \tau \nonumber \\
&&\!+C L(R)\!\!\!\int_0^t \|A^{\beta}_{\alpha}\textrm{e}^{(t-\tau)A_{\alpha}}\|_{L^2(D)}\|u-\bar{u}\|_{L^2(D)}\intd \tau \nonumber \\
&\le& \frac{C L(R)}{\delta}(1-\textrm{e}^{-\delta T})\|u-\bar{u}\|_{H^{\frac{\alpha}{2}}(D)}
+\frac{C L(R)}{1-\beta}T^{1-\beta}\|u-\bar{u}\|_{H^{\frac{\alpha}{2}}(D)},
\end{eqnarray}
where $L(R)$ denote the Lipschitz constant and $\frac12<\beta<1$; now if $T\le T_2$ is choosen small enough,
we get $\|\Phi(u)-\Phi(\bar{u})\|_X<L\|u-\bar{u}\|_X$, $L<1$, making $\Phi$ a contraction mapping from $S$ into itself.
Thus $\Phi$ has a unique fixed point $u$ in $S$, solving ($\ref{solsemi}$). We have thus proved the following result.

\begin{thm}\label{localsol} (Local solution)
The equation ($\ref{semilin}$)
has a unique solution $u\in C^0([0,\,T],\,H_0^{\alpha/2}(D))$, where $T>0$ is appropriately chosen.
\end{thm}

\subsection{Global solution }

Now we examine the global solution   based on the result of
the local existence.
\begin{thm} \label{semithm} (Global solution)
The solution to the nonlocal Ginzberg-Landau equation ($\ref{semilin}$) exists uniquely in the space $C^0(\mathbb{R}^+,\,H_0^{\frac{\alpha}{2}}(D))\cap L_{loc}(\mathbb{R}^+,\,H_0^{\alpha}(D))$. Moreover,
\begin{equation}\label{Haestilocal}
\int_t^{t+1}\|u(x,\,\tau)\|_{H^{\alpha}(D)}^2\intd\tau \le C(|D|,\,\|u_0\|^2_{L^2(D)}).
\end{equation}
\end{thm}
\textbf{Proof}~~It is enough to prove that $\displaystyle{\sup_{0\le t<+\infty}\|u(x,\,t)\|_{H^{\frac{\alpha}{2}}(D)}<+\infty}$, by appropriate energy estimates.

Multiplying $u$  to   both sides of the equation ($\ref{semilin}$)   and integrating, we have
\begin{eqnarray}
\frac12\frac{\intd}{\intd t}\|u\|^2_{L^2(D)}&=&\int_D uA_{\alpha}u\intd x+\int_Du^2 \intd x-\int_Du^4 \intd x \nonumber \\
&\le& -\int_D\int_D (\mathcal{D}^{\ast}(u))^2\intd y \intd x+\int_Du^2 \intd x-\int_Du^4 \intd x \nonumber \\
&=& -\int_D\int_D (\mathcal{D}^{\ast}(u))^2\intd y \intd x+ \frac14|D|-\int_D(u^2-\frac12)^2 \intd x \nonumber \\
&\le& -\int_D\int_D (\mathcal{D}^{\ast}(u))^2\intd y \intd x+\frac14|D|.
\end{eqnarray}
Using the nonlocal Poincar\'e  inequality $\|u\|^2_{L^2(D)}\le C\int_D\int_D (\mathcal{D}^{\ast}(u))^2\intd y \intd x$, we get
\begin{equation}\label{preL2estsemi1}
\frac12\frac{\intd}{\intd t}\|u\|^2_{L^2(D)}+\frac{1}{2C}\|u\|^2_{L^2(D)}+\frac12\int_D\int_D (\mathcal{D}^{\ast}(u))^2\intd y \intd x \le \frac14|D|.
\end{equation}
By the Gronwell's inequality, we conclude that
\begin{equation}\label{L2estsemi1}
\|u\|^2_{L^2(D)} \le \frac12|D|+\|u_0\|^2_{L^2(D)}\textrm{e}^{-\frac{1}{C}t}.
\end{equation}
This implies that  $\displaystyle{\sup_{0\le t<+\infty}\|u(x,\,t)\|_{L^2(D)}<+\infty}$. Integrating ($\ref{preL2estsemi1}$) between $t$ and $t+1$, we obtain
\begin{equation*}
\|u(t+1)\|^2_{L^2(D)}-\|u(t)\|^2_{L^2(D)}+\frac{1}{C}\int_t^{t+1}\|u(\tau)\|^2_{L^2(D)} \intd \tau +\int_t^{t+1}\int_D\int_D (\mathcal{D}^{\ast}(u(\tau)))^2\intd y \intd x \intd\tau \le \frac14|D|.
\end{equation*}
By ($\ref{L2estsemi1}$), we obtain
\begin{equation}\label{preHsestsemi1}
\int_t^{t+1}\int_D\int_D (\mathcal{D}^{\ast}(u(\tau)))^2\intd y \intd x \intd\tau \le \frac14|D|+\|u(t)\|^2_{L^2(D)} \le C(|D|,\,\|u_0\|^2_{L^2(D)}).
\end{equation}


Multiplying $-A_{\alpha}u$ to   both sides of the equation   ($\ref{semilin}$)  and integrating on the domain $D$, we have
\begin{equation}
-\int_D u_tA_{\alpha}u=-\int_D A_{\alpha}u A_{\alpha}u\intd x-\int_D uA_{\alpha}u \intd x+\int_D u^3A_{\alpha}u \intd x.
\end{equation}
By the nonlocal Green's formula \cite{DGLZ2011}, we obtain
\begin{equation}
\frac{\intd}{\intd t}\int_D\int_D (\mathcal{D}^{\ast}(u))^2\intd y \intd x
+\int_D (A_{\alpha}u)^2 \intd x=\int_D\int_D (\mathcal{D}^{\ast}(u))^2\intd y \intd x-\int_D\int_D \mathcal{D}^{\ast}(u^3)\cdot \mathcal{D}^{\ast}(u) \intd y \intd x.
\end{equation}
Now we see that
\begin{eqnarray}
&&-\int_D\int_D \mathcal{D}^{\ast}(u^3)\cdot \mathcal{D}^{\ast}(u) \intd y \intd x=
-\int_D\int_D \mathcal{D}^{\ast}(u)((u^3(y+x)-u^3(x))\alpha(x,\,y) \intd y \intd x \nonumber \\
&=&-\int_D\int_D (\mathcal{D}^{\ast}(u))^2(u^2(y+x)+u(y+x)u(x)+u^2(x)) \intd y \intd x \le 0,
\end{eqnarray}
since $u^2(y+x)+u(y+x)u(x)+u^2(x)$ is nonnegative. Hence,
\begin{equation}\label{preH2sestsemi1}
\frac{\intd}{\intd t}\int_D\int_D (\mathcal{D}^{\ast}(u))^2\intd y \intd x+\int_D (A_{\alpha}u)^2 \intd x
\le \int_D\int_D (\mathcal{D}^{\ast}(u))^2\intd y \intd x.
\end{equation}
Using ($\ref{preHsestsemi1}$) and the uniform Gronwell's inequality (\cite{FP1967, Temam1997}), we obtain
\begin{equation}\label{Hsestsemi1}
\int_D\int_D (\mathcal{D}^{\ast}(u))^2\intd y \intd x \le C(|D|,\,\|u_0\|^2_{L^2(D)}).
\end{equation}
This implies that $\displaystyle{\sup_{0\le t<+\infty}\|u(x,\,t)\|_{H^{\frac{\alpha}{2}}(D)}<+\infty}$.

Integrating ($\ref{preH2sestsemi1}$) between $t$ and $t+1$, we conclude that
\begin{eqnarray*}
&&\int_D\int_D (\mathcal{D}^{\ast}(u(t+1)))^2\intd y \intd x-\int_D\int_D (\mathcal{D}^{\ast}(u(t)))^2\intd y \intd x
+\int_t^{t+1}\int_D (A_{\alpha}u)^2 \intd x \intd\tau  \nonumber \\
&\le& \int_t^{t+1}\int_D\int_D (\mathcal{D}^{\ast}(u(\tau)))^2\intd y \intd x.
\end{eqnarray*}
By ($\ref{Hsestsemi1}$), we obtain
\begin{equation}
\int_t^{t+1}\int_D (A_{\alpha}u)^2 \intd x \intd\tau \le C(|D|,\,\|u_0\|^2_{L^2(D)}),
\end{equation}
and it implies the estimate ($\ref{Haestilocal}$). $\Box$

\begin{rem}
For   linear cases,  a stronger spatial regularity result for the solution can be proved, that is $u(x,\,t)\in C^0(\mathbb{R}^+,\,H^{\alpha}(D))$. In fact,
multiply $A_{\alpha}^2 u$ to both sides of the linear equation $u_t=A_{\alpha}u+u+f(x)$,  we have
\begin{eqnarray*}
&&\frac12 \frac{\intd}{\intd t} \|A_{\alpha} u\|_{L^2(D)}^2 = -\|\mathcal{D}^{\ast} A_{\alpha} u\|_{L^2(D)}^2+\|A_{\alpha}u\|_{L^2(D)}^2
-\int \mathcal{D}^{\ast}(f(x)) \mathcal{D}^{\ast} A_{\alpha} u \\
&&\le C(\|\mathcal{D}^{\ast}(f(x))\|_{L^2(D)}^2+\|u_0\|_{L^2(D)}^2).
\end{eqnarray*}
By the uniform Gronwall inequality, the solution is in $C^0(\mathbb{R}^+,\,H^{\alpha}(D))$.
\end{rem}

\begin{rem} \label{semithm2}
The solution to the following nonlocal semi-linear equation
\begin{equation}\label{semilinear2}
\left\{\begin{array}{ll}
u_t=A_{\alpha}u+u+u^2-u^3,& ~~x\in D,~t>0 \\
u|_{D^c}=0, &\\
u(x,0)=u_0(x), &
\end{array}
\right.
\end{equation}
exists uniquely in  $C^0(\mathbb{R}^+,\,H_0^{\frac{\alpha}{2}}(D))\cap L_{loc}(\mathbb{R}^+,\,H_0^{\alpha}(D))$.
\end{rem}
 Indeed, it is enough to prove that $\displaystyle{\sup_{0\le t<+\infty}\|u(x,\,t)\|_{H^{\frac{\alpha}{2}}(D)}<+\infty}$, by energy estimates.

Multiplying $u$ to   both sides of the equation  (\ref{semilinear2})  and integrating, we have
\begin{eqnarray}
\frac12\frac{\intd}{\intd t}\|u\|^2_{L^2(D)}&=&\int_D uA_{\alpha}u\intd x+\int_Du^2 \intd x
+\int_Du^3 \intd x-\int_Du^4 \intd x \nonumber \\
&\le& -\int_D\int_D (\mathcal{D}^{\ast}(u))^2\intd y \intd x+\int_Du^2 \intd x
+\int_Du^3 \intd x-\int_Du^4 \intd x.
\end{eqnarray}
Using the nonlocal Poincar\'e inequality $\|u\|^2_{L^2(D)}\le C\int_D\int_D (\mathcal{D}^{\ast}(u))^2\intd y \intd x$ and the Cauchy-Schwarz inequality, we get
\begin{eqnarray}
\frac12\frac{\intd}{\intd t}\|u\|^2_{L^2(D)}+\frac{1}{C}\|u\|^2_{L^2(D)}
&\le& \frac32\int_Du^2 \intd x-\frac12\int_Du^4 \intd x \nonumber \\
&\le& \frac98|D|-\frac12\int_D(u^2-\frac32)^2 \intd x \le \frac98|D|.
\end{eqnarray}
By the Gronwell's inequality, we have
\begin{equation}\label{L2est}
\|u\|^2_{L^2(D)} \le (\frac98|D|+\|u_0\|^2_{L^2(D)})\textrm{e}^{-\frac{2}{C}t}.
\end{equation}
This implies $\displaystyle{\sup_{0\le t<+\infty}\|u(x,\,t)\|_{L^2(D)}<+\infty}$.

Multiplying $-A_{\alpha}u$ to   both sides of the equation   (\ref{semilinear2})   and integrating on the domain $D$, we have
\begin{equation}
-\int_D u_tA_{\alpha}u=-\int_D A_{\alpha}u A_{\alpha}u\intd x-\int_D uA_{\alpha}u \intd x-\int_D u^2A_{\alpha}u \intd x+\int_D u^3A_{\alpha}u \intd x.
\end{equation}
By the nonlocal Green's formula again, we obtain
\begin{eqnarray}
&&\frac{\intd}{\intd t}\int_D\int_D (\mathcal{D}^{\ast}(u))^2\intd y \intd x
+\int_D (A_{\alpha}u)^2 \intd x \nonumber \\
&=&\int_D\int_D (\mathcal{D}^{\ast}(u))^2\intd y \intd x+\int_D\int_D \mathcal{D}^{\ast}(u^2)\cdot \mathcal{D}^{\ast}(u) \intd y \intd x
-\int_D\int_D \mathcal{D}^{\ast}(u^3)\cdot \mathcal{D}^{\ast}(u) \intd y \intd x.
\end{eqnarray}
Now we see that
\begin{eqnarray}
&&\int_D\int_D \mathcal{D}^{\ast}(u^2)\cdot \mathcal{D}^{\ast}(u) \intd y \intd x
-\int_D\int_D \mathcal{D}^{\ast}(u^3)\cdot \mathcal{D}^{\ast}(u) \intd y \intd x \nonumber \\
&=&\int_D\int_D \mathcal{D}^{\ast}(u)((u^2(y+x)-u^2(x))\alpha(x,\,y) \intd y \intd x \nonumber \\&&-\int_D\int_D \mathcal{D}^{\ast}(u)((u^3(y+x)-u^3(x))\alpha(x,\,y) \intd y \intd x \nonumber \\
&=&\int_D\int_D (\mathcal{D}^{\ast}(u))^2(u(y+x)+u(x)-u^2(y+x)-u(y+x)u(x)-u^2(x)) \intd y \intd x \nonumber \\
&=&\!\!\!\int_D\int_D (\mathcal{D}^{\ast}(u))^2 (1-\frac12(u(y+x)-1)^2\!-\!\frac12(u(x)-1)^2\!-\!\frac12(u(y+x)+u(x))^2) \intd y \intd x \nonumber \\
&\le& \int_D\int_D (\mathcal{D}^{\ast}(u))^2 \intd y \intd x.
\end{eqnarray}
Hence,
\begin{equation}
\frac{\intd}{\intd t}\int_D\int_D (\mathcal{D}^{\ast}(u))^2\intd y \intd x
\le 2\int_D\int_D (\mathcal{D}^{\ast}(u))^2\intd y \intd x.
\end{equation}
Using the uniform Gronwell's inequality,
\begin{equation}
\int_D\int_D (\mathcal{D}^{\ast}(u))^2\intd y \intd x \le C(T).~~\Box
\end{equation}

\begin{rem} \label{semithm3}
If the nonlinear function $F(\xi)$ is Lipschitz continuous locally, and a dissipative condition holds, i.e. $-C|\xi|^{p}-C\leq F(\xi)\xi  \leq - C|\xi|^{p}+C$  where $p\geq2$, for   $\xi\in\mathbb{R}$,  then the solution to the following nonlocal semi-linear equation
\begin{equation}\label{semilinear3}
\left\{\begin{array}{ll}
u_t=A_{\alpha}u+F(u),& ~~x\in D,~t>0 \\
u|_{D^c}=0, &\\
u(x,0)=u_0(x), &
\end{array}
\right.
\end{equation}
exists uniquely in  $C^0(\mathbb{R}^+,\,L^2(D))\cap L_{loc}(\mathbb{R}^+,\,H_0^{\frac{\alpha}{2}}(D))$.
\end{rem}
The proof is similar to   Theorems $\ref{localsol}$ and $\ref{semithm}$.

\section{A nonlocal Fokker-Planck equation}

We consider the following nonlocal Fokker-Planck equation
\begin{equation}\label{FPE}
\left\{\begin{array}{ll}
u_t=A_{\alpha}u-(b(x)u)_x,& ~~x\in D,~t>0, \\
u|_{D^c}=0, &\\
u(x,0)=u_0(x), &
\end{array}
\right.
\end{equation}
where $A_{\alpha}=-(-\triangle)^{\alpha/2}$ is nonlocal Laplacian operator, where $D=(0,1)$
is an interval in $\mathbb{R}^1$, and   $D^c $ is the complement of $D$.
We will prove the existence and uniqueness of the solution to ($\ref{FPE}$) by the method of vanishing viscosity.

\begin{definition}
For $u \in L^2([0,\,T];H^1(D))$, $u_t \in L^2([0,\,T];L^2(D))$ to be a weak solution to   ($\ref{FPE}$) if
\begin{enumerate}
\renewcommand{\labelenumi}{(\roman{enumi})}
\item
  \begin{equation}\label{weaksol}
  \int_D u_t\varphi\intd x+\int_D (b(x)u)_x\varphi\intd x+\int_D\int_D \mathcal{D}^{\ast}(u)\mathcal{D}^{\ast}(\varphi)\intd y \intd x=0
  \end{equation}
  for each $\varphi\in H_0^1(D)$ and a.e. $0\le t\le T$,
  \item $u(x,\,0)=u_0(x)$.
\end{enumerate}
\end{definition}

We will approximate problem ($\ref{FPE}$) by a parabolic initial-boundary value problem
\begin{equation}\label{viscFPE}
\left\{\begin{array}{ll}
u^{\varepsilon}_t=\varepsilon u^{\varepsilon}_{xx}-(b(x)u^{\varepsilon})_x+A_{\alpha}u^{\varepsilon},& ~~\mbox{in}~D\times (0,\,T] \\ 
u^{\varepsilon}|_{D^c}=0, & ~~\mbox{on}~ D^c\times \{t=0\}\\
u^{\varepsilon}(x,0)=u^{\varepsilon}_0(x), &
\end{array}
\right.
\end{equation}
for $0<\varepsilon\le 1$, $u^{\varepsilon}_0(x)=\eta_{\varepsilon}u_0(x)$ with $\eta$ a modifier. The idea is that for each $\varepsilon>0$,
problem ($\ref{viscFPE}$) has a unique solution $u^{\varepsilon}$ as in \cite[Theorem 3.6, p. 94]{GM1992}. We try to show that as $\varepsilon\to 0$,
$u^{\varepsilon}$ converge to a limit function $u$, which is a weak solution to ($\ref{FPE}$).

\begin{lem}\label{Eele1}(Energy estimates)
If     $b'(x)$ and $b''(x)$ are bounded, then there exists a constant $C$, depending only on domain and the coefficients, such that
\begin{equation}\label{Eevisc}
\max_{0\le t\le T}(\|u^{\varepsilon}(t)\|_{H^1(D)}+\|u^{\varepsilon}_t(t)\|_{L^2(D)})\le C(|D|,\,\|u_0\|^2_{H^1(D)}).
\end{equation}
\end{lem}
\textbf{Proof}~~Multiplying $u^{\varepsilon}$ to   both sides of the equation in ($\ref{viscFPE}$) and integrating on the domain $D$, and using the nonlocal Green's formula,
we have
\begin{eqnarray}
&&\frac12\frac{\intd}{\intd t}\|u^{\varepsilon}\|^2_{L^2(D)}= \varepsilon\int_Du^{\varepsilon}u^{\varepsilon}_{xx}\intd x-\int_D u^{\varepsilon}(b(x)u^{\varepsilon})_x\intd x+\int_D u^{\varepsilon}A_{\alpha}u^{\varepsilon}\intd x \nonumber \\
&=& -\varepsilon\int_D (u_x^{\varepsilon})^2\intd x-\int_D\int_D (\mathcal{D}^{\ast}(u^{\varepsilon}))^2\intd y \intd x-\int_D b(x)u^{\varepsilon}u^{\varepsilon}_x\intd x-\int_Db'(x)(u^{\varepsilon})^2 \intd x \nonumber \\
&=& -\varepsilon\int_D (u_x^{\varepsilon})^2\intd x-\int_D\int_D (\mathcal{D}^{\ast}(u^{\varepsilon}))^2\intd y \intd x-\frac12\int_D b(x)((u^{\varepsilon})^2)_x\intd x-\int_Db'(x)(u^{\varepsilon})^2 \intd x \nonumber \\
&=& -\varepsilon\int_D (u_x^{\varepsilon})^2\intd x-\int_D\int_D (\mathcal{D}^{\ast}(u^{\varepsilon}))^2\intd y \intd x-\frac12\int_Db'(x)(u^{\varepsilon})^2 \intd x.
\end{eqnarray}
Since $b'(x)$ is bounded, i.e. there exists a constant $C$ such that $|b'(x)|<C$ and
\begin{equation}
\frac12\frac{\intd}{\intd t}\|u^{\varepsilon}\|^2_{L^2(D)}\le C\|u^{\varepsilon}\|^2_{L^2(D)}.
\end{equation}
By the Gronwell's inequality, we get
\begin{equation}\label{CLT1L2estsemi}
\max_{0\le t\le T}\|u^{\varepsilon}\|^2_{L^2(D)} \le C\|u_0\|^2_{L^2(D)},
\end{equation}
since $\|u^{\varepsilon}_0\|_{L^2(D)}\le \|u_0\|_{L^2(D)}$.

Multiplying $-u^{\varepsilon}_{xx}$ to   both sides of the equation in ($\ref{viscFPE}$) and integrating on the domain $D$,
we have
\begin{equation}\label{preH1estsemi1}
\frac12\frac{\intd}{\intd t}\|u^{\varepsilon}_x\|^2_{L^2(D)}= -\varepsilon\int_D(u^{\varepsilon}_{xx})^2\intd x+\int_D u^{\varepsilon}_{xx}(b(x)u^{\varepsilon})_x\intd x-\int_D u^{\varepsilon}_{xx}A_{\alpha}u^{\varepsilon}\intd x.
\end{equation}

Suppose $v\in C_0^{\infty}(D)$. Then by the nonlocal Green's formula
\begin{equation}
-\int_D v_{xx}A_{\alpha}v \intd x=\int_D v_x A_{\alpha}v_x=-\int_D\int_D (\mathcal{D}^{\ast}(v_x))^2\intd y \intd x\le 0.
\end{equation}
Then
\begin{eqnarray}
&&\int_D v_{xx}(b(x)v)_x\intd x=\int_D v_{xx}(b(x)v_x) \intd x+\int_D v_{xx}(b'(x)v) \intd x \nonumber \\
&=&-\frac12\int_D b'(x)(v_x)^2 \intd x-\int_D (b'(x)v)_x v_x \intd x \nonumber \\
&=&-\frac32\int_D b'(x)(v_x)^2 \intd x-\int_D b''(x)vv_x \intd x.
\end{eqnarray}
As $b'(x)$ and $b''(x)$ are bounded, we deduce using H\"older inequality that
\begin{equation}
\left|\int_D v_{xx}(b(x)v)_x\intd x\right| \le C(\|v_x\|_{L^2(D)}^2+\|v\|_{L^2(D)}^2).
\end{equation}
Moreover,
\begin{eqnarray}
&&-\int_D u^{\varepsilon}_{xx}A_{\alpha}u^{\varepsilon} \intd x \le 0, \\
&&\left|\int_D u^{\varepsilon}_{xx}(b(x)u^{\varepsilon})_x\intd x\right| \le C(\|u^{\varepsilon}_x\|_{L^2(D)}^2+\|u^{\varepsilon}\|_{L^2(D)}^2).
\end{eqnarray}
Utilizing the above estimates in (\ref{preH1estsemi1}), we obtain
\begin{equation}
\frac{\intd}{\intd t}\|u^{\varepsilon}_x\|^2_{L^2(D)} \le C(\|u^{\varepsilon}_x\|_{L^2(D)}^2+\|u^{\varepsilon}\|_{L^2(D)}^2).
\end{equation}
We next apply Gronwall's inequality, to deduce
\begin{equation}
\max_{0\le t\le T}\|u^{\varepsilon}_x\|^2_{L^2(D)} \le C(\|u_0\|^2_{L^2(D)}+\|(u_0)_x\|_{L^2(D)}^2),
\end{equation}
using $\|(u^{\varepsilon}_0)_x\|_{L^2(D)}\le \|(u_0)_x\|_{L^2(D)}$ and the $L^2$-estimate ($\ref{CLT1L2estsemi}$).

Differentiating with respect to $t$, multiplying $u^{\varepsilon}_t$ to   both sides of the equation in ($\ref{viscFPE}$) and integrating on the domain $D$, we have
\begin{equation}\label{preH1estsemi1t}
\frac12\frac{\intd}{\intd t}\int_D (u^{\varepsilon}_t)^2 \intd x =\varepsilon \int_D u^{\varepsilon}_t(u^{\varepsilon}_t)_{xx} \intd x -\int_D u^{\varepsilon}_t(b(x)u^{\varepsilon})_{tx} \intd x+\int_D u^{\varepsilon}_tA_{\alpha}u^{\varepsilon}_t \intd x.
\end{equation}

Suppose $v\in C_0^{\infty}(\overline{D\times(0,\,T)})$.   By the classical and nonlocal Green's formulas
\begin{eqnarray*}
&&\varepsilon \int_D v_t(v_t)_{xx} \intd x=-\varepsilon \int_D ((v_t)_x)^2 \intd x\le 0, \\
&&\int_D v_tA_{\alpha}v_t \intd x-\int_D\int_D (\mathcal{D}^{\ast}(v_t))^2\intd y \intd x\le 0.
\end{eqnarray*}
Thus
\begin{equation*}
-\int_D v_t(b(x)v)_{tx} \intd x=-\int_D b'(x)(v_t)^2 \intd x-\int_D v_tb(x)(v_t)_x \intd x=-\frac12\int_D b'(x)(v_t)^2 \intd x.
\end{equation*}
As $b'(x)$ is bounded, we deduce
\begin{equation*}
\left|\int_D v_t(b(x)v)_{tx} \intd x\right| \le C\|v_t\|_{L^2(D)}^2.
\end{equation*}
Also note that
\begin{eqnarray*}
&&\varepsilon \int_D u^{\varepsilon}_t(u^{\varepsilon}_t)_{xx} \intd x \le 0,~~~\int_D u^{\varepsilon}_tA_{\alpha}u^{\varepsilon}_t \intd x \le 0, \\
&&\left|\int_D u^{\varepsilon}_t(b(x)u^{\varepsilon})_{tx} \intd x\right| \le C\|u^{\varepsilon}_t\|_{L^2(D)}^2.
\end{eqnarray*}
Utilizing the above estimates   in (\ref{preH1estsemi1t}), we obtain
\begin{equation}
\frac{\intd}{\intd t}\|u^{\varepsilon}_t\|^2_{L^2(D)} \le C\|u^{\varepsilon}_t\|_{L^2(D)}^2.
\end{equation}
We finally apply the Gronwall's inequality, to deduce
\begin{equation}
\max_{0\le t\le T}\|u^{\varepsilon}_t\|^2_{L^2(D)} \le C\|(u_0)_t\|_{L^2(D)},
\end{equation}
where we have used the fact $\|(u^{\varepsilon}_0)_t\|_{L^2(D)}\le \|(u_0)_t\|_{L^2(D)}$. $\Box$

\begin{thm}\label{exisFPE} (Weak solution)
There exist a unique weak solution to the initial-boundary problem ($\ref{FPE}$).
\end{thm}
\textbf{Proof}~~We first prove the existence of the weak solution. According to the energy estimate ($\ref{Eele1}$), there exists a subsequence $\varepsilon_k \to 0$ and a function $u\in L^2((0,\,T),\,H^1(D))$, such that $u_t\in L^2((0,\,T),\,L^2(D))$, with
\begin{equation}
\left\{\begin{array}{ll}
u^{\varepsilon_k}\rightharpoonup u & ~~\textrm{in}~L^2((0,\,T),\,H^1(D)), \\
u^{\varepsilon_k}_t\rightharpoonup u_t & ~~\textrm{in}~L^2((0,\,T),\,L^2(D)).
\end{array}\right.
\end{equation}

Choose a function $\varphi\in C^1([0,\,T];\,H^1(D))$. Then from equation ($\ref{viscFPE}$) we deduce
\begin{equation}\label{beflimit}
\int_0^T\left(\int_D u^{\varepsilon}_t\varphi \intd x +\varepsilon\int_Du^{\varepsilon}_{x}\varphi_x\intd x+\int_D (b(x)u^{\varepsilon})_x\varphi\intd x+\int_D \mathcal{D}^{\ast}(u^{\varepsilon})\mathcal{D}^{\ast}(\varphi) \intd x\right) \intd t=0.
\end{equation}
Letting $\varepsilon=\varepsilon_k \to 0$, we get
\begin{equation}\label{afterlimited}
\int_0^T\left(\int_D u_t\varphi \intd x +\int_D (b(x)u)_x\varphi\intd x+\int_D \mathcal{D}^{\ast}(u)\mathcal{D}^{\ast}(\varphi) \intd x\right) \intd t=0.
\end{equation}
The identity above is also valid for all $\varphi\in C([0,\,T];\,H^1(D))$, and so
\begin{equation*}
\int_D u_t\varphi \intd x +\int_D (b(x)u)_x\varphi\intd x+\int_D \mathcal{D}^{\ast}(u)\mathcal{D}^{\ast}(\varphi) \intd x=0
\end{equation*}
for a.e.   $t$ and each $\varphi \in H^1(D)$.

Now assume $\varphi(T)=0$. Then ($\ref{beflimit}$) implies
\begin{eqnarray*}
&&\int_0^T\left(-\int_D u^{\varepsilon}\varphi_t \intd x +\varepsilon\int_Du^{\varepsilon}_{x}\varphi_x\intd x+\int_D (b(x)u^{\varepsilon})_x\varphi\intd x+\int_D \mathcal{D}^{\ast}(u^{\varepsilon})\mathcal{D}^{\ast}(\varphi)
\intd x\right) \intd t \\
&=&\int_D u^{\varepsilon}_0\varphi(0) \intd x.
\end{eqnarray*}
Letting  $\varepsilon=\varepsilon_k \to 0$, we obtain
\begin{equation*}
\int_0^T\left(-\int_D u\varphi \intd x +\int_D (b(x)u)_x\varphi\intd x+\int_D \mathcal{D}^{\ast}(u)\mathcal{D}^{\ast}(\varphi)
\intd x\right) \intd t=\int_D u_0\varphi(0) \intd x.
\end{equation*}
Integrating by parts in ($\ref{afterlimited}$) gives us the identity
\begin{equation*}
\int_0^T\left(-\int_D u_t\varphi \intd x +\int_D (b(x)u)_x\varphi\intd x+\int_D \mathcal{D}^{\ast}(u)\mathcal{D}^{\ast}(\varphi) \intd x\right) \intd t=\int_D u(x,0)\varphi(0) \intd x.
\end{equation*}
Consequently $u(x,0)=u_0(x)$, as $\varphi(0)$ is arbitrary.

Now it suffices to demonstrate the only weak solution of ($\ref{FPE}$) with $u_0 \equiv 0$ is $u \equiv 0$. To verify this, note that
\begin{equation}\label{weaksol1}
\int_D u_tu \intd x +\int_D (b(x)u)_xu\intd x+\int_D \mathcal{D}^{\ast}(u)\mathcal{D}^{\ast}(u) \intd x=0~~ \mbox{for~a.~e.}~
0\le t\le T.
\end{equation}
Since $b'(x)$ is bounded, we  calculate from ($\ref{weaksol1}$) that
\begin{equation}
\frac12\frac{\intd}{\intd t}\|u\|^2_{L^2(D)}\le C\|u\|^2_{L^2(D)}.
\end{equation}
Hence the Gronwall's inequality forces $\|u\|^2_{L^2(D)}=0 $ for $0\le t\le T$, as $u(x,\,0)=0$.  This implies uniqueness. $\Box$

\appendix
\section{The Uniform Gronwall's inequality}

We present the uniform Gronwall's inequality that was implicitly used for the first time by   Foias and   Prodi \cite{FP1967} in the context of the Navier-Stokes equations.

Let $g$, $h$, $y$ be three locally integrable functions on $(t_0,\,+\infty)$ that satisfy
\begin{equation}\label{Gronwallcondition}
\frac{\intd y}{\intd t} \le gy+h~~~\mbox{for}~~t \ge t_0,
\end{equation}
the function $\intd y/\intd t$ being also locally integrable. For the usual Gronwall's inequality we multiply ($\ref{UGIcondition}$) by
$$
\textrm{exp}\left(-\int_{t_0}^t g(\tau) \intd \tau \right),
$$
and observe that the resulting inequality reads
$$
\frac{\intd}{\intd t}\left(y(t)\textrm{exp}\left(-\int_{t_0}^t g(\tau) \intd \tau \right)\right)
\le h(t)\textrm{exp}\left(-\int_{t_0}^t g(\tau) \intd \tau \right).
$$
Hence, by integration between $t_0$ and $t$,
\begin{equation}
y(t) \le y(t_0)\textrm{exp}\left(\int_{t_0}^t g(\tau) \intd \tau \right)+ \int_{t_0}^t h(s)
\textrm{exp}\left(-\int_{t_0}^t g(\tau) \intd \tau \right) \intd s, ~~~~t \ge t_0.
\end{equation}
This is the usual Gronwall's inequality which is useful for bounded values of $t$. When $t\to\infty$,
this relation is not sufficient for our purposes since it allows an exponential growth of $y$; for instant, for $y \ge 0$, $h=0$, $g=1$, we find
$$
y(t) \le y(t_0)\textrm{exp}(t-t_0).
$$

We now present an alternative form of this inequality that provides (under slightly stronger assumptions) a bound valid uniformly for $t \ge t_0$.
\begin{lem} (Uniform Gronwall's inequality)
Let $g$, $h$, $y$ be three positive locally integrable functions on $(t_0,\,+\infty)$ such that $y'$ is locally integrable on $(t_0,\,+\infty)$, and
\begin{equation}\label{UGIcondition}
\frac{\intd y}{\intd t} \le gy+h~~~\mbox{for}~~t \ge t_0,
\end{equation}
\begin{equation}
\int_t^{t+r} g(s) \intd s \le a_1,~\int_t^{t+r} h(s) \intd s \le a_2,~\int_t^{t+r} y(s) \intd s \le a_3,~~~\mbox{for}~~t \ge t_0,
\end{equation}
where $r$, $a_1$, $a_2$, $a_3$ are positive constants. Then
\begin{equation}\label{UGIresult}
y(t+r) \le \left(\frac{a_3}{r}+a_2\right) \exp(a_1),~~~\forall\,t\ge t_0.
\end{equation}
\end{lem}





\end{document}